\documentclass{amsart}

\newtheorem{theorem}{Theorem}[section]
\newtheorem{lemma}[theorem]{Lemma}

\newtheorem{corollary}[theorem]{Corollary}
\newtheorem{question}[theorem]{Question}
\newtheorem{conjecture}[theorem]{Conjecture}

\theoremstyle{definition}

\newtheorem{example}[theorem]{Example}
\newtheorem{algorithm}[theorem]{Algorithm}

\theoremstyle{remark}
\newtheorem{remark}[theorem]{Remark}

\numberwithin{equation}{section}

\begin{document}

\title[Cyclic Resultants]{Cyclic Resultants}

\author{Christopher J. Hillar}
\address{Department of Mathematics, University of California, Berkeley, CA 94720.}
\email{chillar@math.berkeley.edu}
\thanks{This work is supported under a National
Science Foundation Graduate Research Fellowship.}



\begin{abstract}
We characterize polynomials having the same set of nonzero cyclic
resultants.  Generically, for a polynomial $f$ of degree $d$,
there are exactly $2^{d-1}$ distinct degree $d$ polynomials with
the same set of cyclic resultants as $f$.  However, in the generic
monic case, degree $d$ polynomials are uniquely determined by
their cyclic resultants. Moreover, two reciprocal
(``palindromic'') polynomials giving rise to the same set of
nonzero cyclic resultants are equal.  In the process, we also prove a unique factorization
result in semigroup algebras involving products of binomials.  Finally, we discuss how our results yield algorithms for explicit reconstruction of polynomials from their cyclic resultants.
\end{abstract} \maketitle

\section{Introduction}

The $m$-th cyclic resultant of a univariate polynomial $f \in \mathbb C[x]$ is \[r_m = \text{Res}(f,x^m-1).\] We are primarily interested here in the fibers of the map $r: \mathbb C[x] \to
\mathbb C^{\mathbb N}$ given by $f \mapsto
\left(r_m\right)_{m=0}^{\infty}$.  In particular, what are the
conditions for two polynomials to give rise to the same set of
cyclic resultants?  For technical reasons, we will only consider
polynomials $f$ that do not have a root of unity as a zero. With
this restriction, a polynomial will map to a set of all nonzero
cyclic resultants.  Our main result gives a complete answer to this question.  

\begin{theorem}\label{characterizationthmchar0}
Let f and g be polynomials in $\mathbb C[x]$. Then, $f$ and $g$ generate the same
sequence of nonzero cyclic resultants if and only if there exist $u, v 
\in \mathbb C[x]$ with
$u(0)\not =0$ and nonnegative integers $l_1,l_2$ such that 
$\deg(u) \equiv l_2-l_1 \ (\text{\rm mod} \ 2)$, and 
\begin{equation*}
\begin{split}
f(x) =&  \ (-1)^{l_2-l_1}x^{l_1}v(x)u(x^{-1})x^{\deg(u)} \\
g(x) =&  \ x^{l_2}v(x)u(x). \\
\end{split}
\end{equation*}
\end{theorem}

\begin{remark}
All our results involving $\mathbb C$ hold over any algebraically closed field of characteristic zero.
\end{remark}

Although the theorem statement appears somewhat technical, we
present a natural interpretation of the result.  Suppose that
$g(x) = x^{l_2}v(x)u(x)$ is a factorization as above of a polynomial $g$
with nonzero cyclic resultants.  Then, another polynomial $f$
giving rise to this same sequence of resultants is obtained from
$v$ by multiplication with the reversal $u(x^{-1})x^{\deg(u)}$ of $u$ and a factor
$(-1)^{\deg(u)}x^{l_1}$ in which $l_1 \equiv l_2 - \deg(u)  \ (\text{\rm mod} \ 2)$. 
In other words, $f(x) = (-1)^{\deg(u)}x^{l_1}v(x)u(x^{-1})x^{\text{deg}(u)}$, 
and all such $f$ must
arise in this manner.

\begin{example}
One can check that the polynomials 
\begin{equation*}
\begin{split}
f(x) = 
& \ {x}^{3}-10\,{x}^{2}+31\,x-30 \\
g(x) = & \ 15\,{x}^{5}-38\,{x}^{4}+17\,{x}^{3}-2\,{x}^{2} \\
\end{split}
\end{equation*} both generate
the same cyclic resultants.  This follows from the factorizations
\begin{equation*}
\begin{split}
f(x) =& \ (x-2)\left(15x^2-8x+1 \right) \\
g(x) =&  \ x^2(x-2)\left(x^2-8x+15\right). \qed  \\
\end{split}
\end{equation*}
\end{example}

One motivation for the study of cyclic resultants comes from the
theory of dynamical systems.  Sequences of the form $r_m$ arise as
the cardinalities of sets of periodic points for toral
endomorphisms. Let $A$ be a $d$-by-$d$ integer matrix and let $X = \mathbb T^{d} = \mathbb R^{d}/\mathbb Z^{d}$ denote the $d$-dimensional additive torus. Then, the matrix $A$ acts on $X$ by multiplication mod $1$; that is, it defines a map $T: X \to X$ given by
\[T(\mathbf{x}) = A\mathbf{x} \mod \mathbb Z^d.\]  Let $\text{Per}_m(T) =
\{\mathbf{x} \in \mathbb T^{d} : T^m(\mathbf{x}) = \mathbf{x}\}$
be the set of points fixed under the map $T^m$. Under the
ergodicity condition that no eigenvalue of $A$ is a root of unity, it
follows (see \cite{Ward}) that \[|r_m(f)| = |\text{Per}_m(T)| =
|\det(A^m-I)|,\] in which $I$ is the $d$-by-$d$ identity matrix, and $f$ is the characteristic polynomial of $A$.  As a consequence of our results, we characterize when the sequence
$|\text{Per}_m(T)|$ determines the spectrum of the linear map $A$ lifting $T$ (see Corollary \ref{toralcor}).

In connection with number theory, cyclic resultants were also studied
by Pierce and Lehmer \cite{Ward} in the hope of using them to
produce large primes.  As a simple example, the Mersenne 
numbers $M_m = 2^m - 1$ arise as cyclic resultants of the polynomial 
$f(x) =x-2$.  Indeed, the map $T(x) = 2x \mod 1$ has precisely $M_m$ 
points of period $m$.  Further motivation comes from knot theory \cite{knots},
Lagrangian mechanics \cite{guillemin,Zworski}, and, more recently,
in the study of amoebas of varieties \cite{purbhoo} and quantum computing \cite{kedlaya}.

The principal result in the direction of our main characterization
theorem was discovered by Fried \cite{fried} although certain
implications of Fried's result were known to Stark
\cite{duistermaat}.  Our approach is a refinement and generalization of the one found in \cite{fried}. 
Given a polynomial $f = a_0x^d+a_1x^{d-1}+\cdots+a_d$ of
degree $d$, the \textit{reversal} of $f$ is the polynomial
$x^df(1/x)$. Additionally, $f$ is called \textit{reciprocal} if
$a_i = a_{d-i}$ for $0 \leq i \leq d$ (sometimes such a polynomial
is called \textit{palindromic}). Alternatively, $f$ is reciprocal
if it is equal to its own reversal.  Fried's result may be stated
as follows.  It will be a corollary of Theorem \ref{charthmrealcase} below (the real version of Theorem \ref{characterizationthmchar0}).

\begin{corollary}[Fried]\label{friedtheorem}
Let $p(x) = a_0 x^d + \cdots + a_{d-1} x + a_d \in \mathbb R[x]$
be a real reciprocal polynomial of even degree $d$ with $a_0 > 0$,
and let $r_m$ be the $m$-th cyclic resultants of $p$. Then,
$|r_m|$ uniquely determine this polynomial of degree $d$ as long
as the $r_m$ are never $0$.
\end{corollary}

The following is a direct corollary of our main theorem to the
generic case.

\begin{corollary}\label{genericcor}
Let $g$ be a generic
polynomial in $\mathbb C[x]$ of degree $d$.  Then, there are exactly
$2^{d-1}$ degree $d$ polynomials with the same set of
cyclic resultants as $g$.
\end{corollary}

\begin{proof}
If $g$ is generic, then $g$ will not have a root of unity as a
zero nor will $g(0) = 0$.  Theorem \ref{characterizationthmchar0},
therefore, implies that any other degree $d$ polynomial $f \in
\mathbb C[x]$ giving rise to the same set of cyclic resultants
is determined by choosing an even cardinality subset of the roots
of $g$.  Such polynomials will be distinct since $g$ is generic.
Since there are $2^d$ subsets of the roots of $g$ and half of them
have even cardinality, the theorem follows.
\end{proof}

\begin{example}
Let $g(x) = (x-2)(x-3)(x-5) = {x}^{3}-10\,{x}^{2}+31\,x-30$. Then,
there are $2^{3-1}-1 = 3$ other degree $3$ polynomials with the
same set of cyclic resultants as $g$.  They are:
\[15\,{x}^{3}-38\,{x}^{2}+17\,x-2\] \[10\,{x}^{3}-37\,{x}^{2}+22\,x-3\]
\[6\,{x}^{3}-35\,{x}^{2}+26\,x-5.\qed \]
\end{example}

If one is interested in the case of generic monic polynomials,
then Theorem \ref{characterizationthmchar0} also implies the
following uniqueness result.

\begin{corollary}\label{genericcormonic}
The set of cyclic resultants determines $g$ for generic monic $g \in \mathbb C[x]$ of degree $d$.
\end{corollary}

\begin{proof}
Again, since $g$ is generic, it will not have a root of unity as a
zero nor will $g(0) = 0$.  Theorem \ref{characterizationthmchar0}
forces a constraint on the roots of $g$ for there to be a
different monic polynomial $f$ with the same set of cyclic resultants as
$g$.  Namely, a subset of the roots of $g$ has product $1$, a
non-generic situation.
\end{proof}

As to be expected, there are analogs of Theorem
\ref{characterizationthmchar0} and Corollary \ref{genericcormonic}
to the real case involving absolute values.

\begin{theorem}\label{charthmrealcase}
Let $f$ and $g$ be polynomials in $\mathbb R[x]$. If $f$ and $g$
generate the same sequence of nonzero cyclic resultant absolute
values, then there exist $u,v \in \mathbb C[x]$ with $u(0) \neq 0$
and nonnegative integers $l_1, l_2$ such that
\begin{equation*}
\begin{split}
f(x) =&  \ \pm x^{l_1}v(x)u(x^{-1})x^{\deg(u)} \\
g(x) =&  \ x^{l_2}v(x)u(x). \\
\end{split}
\end{equation*}
\end{theorem}

\begin{corollary}\label{genericcorreal}
The set of cyclic resultant absolute values determines $g$ for generic monic $g \in \mathbb R[x]$ of degree $d$.
\end{corollary}

The generic real case without the monic assumption is more subtle than that of Corollary \ref{genericcor}. The
difficulty is that we are restricted to polynomials in $\mathbb
R[x]$. However, there is the following

\begin{corollary}\label{genericrealthm}
Let $g$ be a generic polynomial in the set of degree $d$ elements of $\mathbb R[x]$ with at most one real root. Then there are exactly $2^{\lceil d/2 \rceil+1}$ degree
$d$ polynomials in $\mathbb R[x]$ with the same set of cyclic
resultant absolute values as $g$.
\end{corollary}

\begin{proof}
If $d$ is even, then the hypothesis implies that all of the roots of
$g$ are nonreal.  In particular, it follows from Theorem
\ref{charthmrealcase} (and genericity) that any other degree $d$
polynomial $f \in \mathbb R[x]$ giving rise to the same set of
cyclic resultant absolute values is determined by choosing a
subset of the $d/2$ pairs of conjugate roots of $g$ and a sign.
This gives us a count of $2^{d/2+1}$ distinct real polynomials.
When $d$ is odd, $g$ has exactly one real root, and a
similar counting argument gives us $2^{\lceil d/2 \rceil+1}$ for
the number of distinct real polynomials in this case.  This proves
the corollary.
\end{proof}

A surprising consequence of this result is that the number of
polynomials with equal sets of cyclic resultant absolute values can be significantly smaller than the number predicted by Corollary \ref{genericcor}.

\begin{example}
Let $g(x) = (x-2)(x+i+2)(x-i+2) = {x}^{3}+2\,{x}^{2}-3\,x-10$.
Then, there are $2^{\lceil 3/2 \rceil+1}-1 = 7$ other degree $3$
real polynomials with the same set of cyclic resultant absolute
values as $g$.  They are: \[-{x}^{3}-2\,{x}^{2}+3\,x+10, \ \pm(-2\,{x}^{3}-7\,{x}^{2}-6\,x+5),\]
\[\pm(5\,{x}^{3}-6\,{x}^{2}-7\,x-2), \  \pm(-10\,{x}^{3}-3\,{x}^{2}+2\,x+1).\]  It is important to realize that while
\begin{equation*}
\begin{split}
f(x) = & \ (1-2x)(1+(i+2)x)(x-i+2) \\
= & \ \left( -4-2\,i \right) {x}^{3}- \left( 10-i \right) {x}^{2}+
\left( 2+2\,i \right) x+2-i
\end{split}
\end{equation*}
has the same set of actual cyclic resultants (by Theorem
\ref{characterizationthmchar0}), it does not appear in the count
above since it is not in $\mathbb R[x]$. 
\qed \end{example}

As an illustration of the usefulness of Theorem
\ref{characterizationthmchar0}, we prove a uniqueness result
involving cyclic resultants of reciprocal polynomials. Fried's
result also follows in the same way using Theorem
\ref{charthmrealcase} in place of Theorem
\ref{characterizationthmchar0}.

\begin{corollary}\label{recipcorollary}
Let $f$ and $g$ be reciprocal polynomials with equal sets of
nonzero cyclic resultants.  Then, $f = g$.
\end{corollary}

\begin{proof}
Let $f$ and $g$ be reciprocal polynomials having the same set of
nonzero cyclic resultants. Applying Theorem
\ref{characterizationthmchar0}, it follows that $d =$ deg($f$) =
deg($g$) and that
\begin{equation*}
\begin{split}
f(x) =& \ v(x)u(x^{-1})x^{\text{deg}(u)} \\
g(x) =& \ v(x)u(x) \\
\end{split}
\end{equation*}
($l_1 = l_2 = 0$ since $f(0),g(0) \neq 0$). But then,
\begin{equation*}
\begin{split}
\frac{u(x^{-1})}{u(x)}x^{\text{deg}(u)} &= \frac{f(x)}{g(x)} \\
& = \frac{x^df(x^{-1})}{x^dg(x^{-1})} \\
& = \frac{u(x)}{u(x^{-1})}x^{-\text{deg}(u)}. \\
\end{split}
\end{equation*}
In particular, $u(x) = \pm u(x^{-1})x^{\text{deg}(u)}$.  If $u(x)
= u(x^{-1})x^{\text{deg}(u)}$, then $f = g$ as desired.  In the
other case, it follows that $f = -g$.  But then Res($f$,$x-1$) $=$
Res($g$,$x-1$) $=$ $-$Res($f$,$x-1$) is a contradiction to $f$
having all nonzero cyclic resultants.  This completes the proof.
\end{proof}

We now state the application to toral endomorphims discussed in the introduction.

\begin{corollary}\label{toralcor}
Let $T$ be an ergodic, toral endomorphism induced by a $d$-by-$d$ integer matrix $A$. If there is no subset of the eigenvalues of $A$ with product $\pm 1$, then the sequence $|\text{\rm{Per}}_m(T)|$ determines the spectrum of the linear map that defines $T$.
\end{corollary}

\begin{proof}
Suppose that $T'$ is another toral endomorphism induced by an integral $d$-by-$d$ matrix $B$ such that  \[|\text{Per}_m(T)| =  |\text{Per}_m(T')|.\]  Let $f$ and $g$ be the characteristic polynomials of $A$ and $B$, respectively.  From the hypothesis of the corollary and the statement of Theorem \ref{charthmrealcase}, it follows that $f$ and $g$ must be equal.  In particular, the eigenvalues of the matrices $A$ and $B$ coincide, completing the proof.
\end{proof}

\begin{remark}
We note that a more complete characterization is possible using the results of Theorem \ref{charthmrealcase}, however, the statement is more technical and not very enlightening.
\end{remark}

When a degree $d$ polynomial is uniquely determined by its sequence of cyclic resultants, it is natural to ask for an algorithm that performs the reconstruction.  In several applications, moreover, explicit inversion using small numbers of resultants is desired (see, for instance, \cite{Zworski, kedlaya}).  In Section \ref{reconstruct}, we describe a method that inverts the map $r$ using the first $2^{d+1}$ cyclic resultants.  Empirically, however, only $d+1$ resultants suffice, and a conjecture by Sturmfels and Zworski would imply that this is always the case.  As evidence for this conjecture, we provide explicit reconstructions for several small examples.

The rest of the paper is organized as follows.  In Section \ref{binfact}, we make a digression into the theory of semigroup algebras and binomial factorizations.  The unique factorization result discussed there (Theorem \ref{mainfactortheorem}) will form a crucial component in proving Theorem \ref{characterizationthmchar0}.  The subsequent chapter deals with algebraic properties of cyclic resultants, and Section \ref{reconstruct} concludes with proofs of our main cyclic resultant characterization theorems.  Finally, in the last section, we discuss algorithms for reconstruction.

\section{Binomial Factorizations}\label{binfact}

We now switch to the seemingly unrelated topic of binomial
factorizations in semigroup algebras.  The relationship to cyclic
resultants will become clear later.  Let $A$ be a finitely
generated abelian group and let $a_1,\ldots,a_n$ be distinguished
generators of $A$. Let $Q$ be the semigroup generated by
$a_1,\ldots,a_n$. The \textit{semigroup
algebra} $\mathbb C[Q]$ is the $\mathbb C$-algebra with vector space basis
$\{\textbf{s}^a : a \in Q\}$ and multiplication defined by
$\textbf{s}^a \cdot \textbf{s}^b = \textbf{s}^{a+b}$. Let $L$
denote the kernel of the homomorphism $\mathbb Z^n$ onto $A$. The
\textit{lattice ideal} associated with $L$ is the following ideal
in $S = \mathbb C[x_1,\ldots,x_n]$: \[I_L = \langle x^u-x^v \ : \ u,v \in
\mathbb N^n \text{ with } u-v \in L \rangle.\]

It is well-known that $\mathbb C[Q] \cong S/I_L$ (e.g. see
\cite{CCA}).  We are primarily concerned here with certain kinds
of factorizations in $\mathbb C[Q]$.

\begin{question}
When is a product of binomials in $\mathbb C[Q]$ equal to another product
of binomials?
\end{question}

The answer to this question turns out to be fundamental for the
study of cyclic resultants.  Our main result in this direction is
a certain kind of unique factorization of binomials in $\mathbb C[Q]$.

\begin{theorem}\label{mainfactortheorem}
Let $\alpha \in \mathbb C$ and suppose that \[ \textbf{s}^a \prod\limits_{i = 1}^e {\left(
{\textbf{s}^{u_i } - \textbf{s}^{v_i } } \right) = \alpha
\textbf{s}^b } \prod\limits_{i = 1}^f {\left( {\textbf{s}^{x_i } -
\textbf{s}^{y_i } } \right)}\] are two factorizations of binomials
in the ring $\mathbb C[Q]$. Furthermore, suppose that for each $i$, the difference $u_i -
v_i$ ($\text{resp. } x_i - y_i$) has infinite order as an element of $A$. Then,
$\alpha = \pm1$, $e = f$, and up to permutation, for each $i$,
there are elements $c_i,d_i \in Q$ such that
$\textbf{s}^{c_i}(\textbf{s}^{u_i}-\textbf{s}^{v_i}) = \pm
\textbf{s}^{d_i}(\textbf{s}^{x_i}-\textbf{s}^{y_i})$.
\end{theorem}

Of course, when each side has a factor of zero, the theorem fails.
There are other obstructions, however, that make necessary the
supplemental hypotheses concerning order.  For example, when $A = \mathbb Z/2\mathbb Z$, we have $\mathbb C[Q] = \mathbb C[A] \cong \mathbb Q[s]/ \langle s^2-1 \rangle$, and it is easily verified that  \[(1-s)(1-s) = 2(1-s).\]

One might also wonder what happens when the binomials are not of the
form $\textbf{s}^u-\textbf{s}^v$.  The following example exhibits
some of the difficulty in formulating a general statement.

\begin{example}\label{nonexample}
$L = \{(0,b) \in \mathbb Z^2 : b \text{ is even}\}$, $I_L =
\langle s^2-1 \rangle \subseteq \mathbb C[s,t]$, $A = \mathbb Z \oplus
\mathbb Z/2\mathbb Z$, $Q = \mathbb N \oplus \mathbb Z/2\mathbb
Z$.  Then, \[(1-t^4) =
(1-st)(1+st)(1-ist)(1+ist) = (1-st^2)(1+st^2)\] are three
different binomial factorizations of the same semigroup algebra
element.
\qed \end{example}

We now are in a position to outline our strategy for
characterizing those polynomials $f$ and $g$ having the same set
of nonzero cyclic resultants (this strategy is similar to the one
employed in \cite{fried}). Given a polynomial $f$ and its sequence
of $r_m$, we construct the generating function $E_f(z) =
\exp{\left(-\sum_{m \geq 1} {r_m \frac{z^m}{m}}\right)}$. This
series turns out to be rational with coefficients depending
explicitly on the roots of $f$. Since $f$ and $g$ are assumed to
have the same set of $r_m$, it follows that their corresponding
rational functions $E_f$ and $E_g$ are equal. Let $G$ be the
(multiplicative) group of units of $\mathbb C$.
Then, the divisors of these two rational functions are group ring
elements in $\mathbb Z[G]$, and their equality forces a certain
binomial group ring factorization that is analyzed explicitly. The
main results in the introduction follow from this final analysis.

To prove our factorization result, we will pass to the full group
algebra $\mathbb C[A]$. As above, we represent elements $\tau \in \mathbb C[A]$ as
$\tau = \sum_{i=1}^{m} { \alpha_i \textbf{s}^{g_i} }$, in which
$\alpha_i \in \mathbb C$ and $g_i \in A$. The following lemma is quite
well-known.

\begin{lemma}\label{zerodivlemma}
If $0 \neq \alpha \in \mathbb C$ and $g \in A$ has infinite order, then
$1-\alpha\textbf{s}^g \in \mathbb C[A]$ is not a zero-divisor.
\end{lemma}

\begin{proof}
Let $0 \neq \alpha \in \mathbb C, g \in A$ and $\tau = \sum_{i=1}^{m} {
\alpha_i \textbf{s}^{g_i} } \neq 0$ be such that \[\tau =
\alpha\textbf{s}^g\tau = \alpha^2 \textbf{s}^{2g}\tau =
\alpha^3 \textbf{s}^{3g}\tau  = \cdots.\]  Suppose that $\alpha_1
\neq 0$.  Then, the elements $\textbf{s}^{g_1},
\textbf{s}^{g_1+g}, \textbf{s}^{g_1+2g}, \ldots$ appear in $\tau$
with nonzero coefficient, and since $g$ has infinite order, these
elements are all distinct.  It follows, therefore, that $\tau$
cannot be a finite sum, and this contradiction finishes the proof.
\end{proof}

Since the proof of the main theorem involves multiple steps, we
record several facts that will be useful later.  The first result
is a verification of the factorization theorem for a special case.

\begin{lemma}\label{polyequality1var}
Fix an abelian group $C$.  Let $\mathbb C[C]$ be the group algebra with $\mathbb C$-vector
space basis given by $\{\textbf{s}^c : c \in C\}$ and set $R =
\mathbb C[C][t,t^{-1}]$.  Suppose that $c_i,d_i, b
\in C$, $m_i,n_i$ are nonzero integers, $q
\in \mathbb Z$, and $z \in \mathbb C$ are such that
\[\prod_{i=1}^{e} {(1-\textbf{s}^{c_i}t^{m_i})} =
z\textbf{s}^{b}t^{q} \prod_{i=1}^{f}
{(1-\textbf{s}^{d_i}t^{n_i})}\] holds in R.  Then, $e = f$ and
after a permutation, for each $i$, either $\textbf{s}^{c_i}t^{m_i}
= \textbf{s}^{d_i}t^{n_i}$ or $\textbf{s}^{c_i}t^{m_i}=
\textbf{s}^{-d_i}t^{-n_i}$.
\end{lemma}

\begin{proof}
Let $\text{sgn}: \mathbb Z \setminus \{0\} \to \{-1,1\}$ denote
the standard sign map $\text{sgn}(n) = n/|n|$ and set $\gamma =
z\textbf{s}^{b}t^{q}$. Rewrite the left-hand side of the given
equality as:
\begin{equation*}
\begin{split}
\prod_{i=1}^{e} {(1-\textbf{s}^{c_i}t^{m_i})} &=
\prod_{\text{sgn}(m_i) = -1} {-\textbf{s}^{c_i}t^{m_i}} \
 \prod_{i=1}^{e}{\left(1-\textbf{s}^{\text{sgn}(m_i)c_i}t^{|m_i|}\right)}.
\end{split}
\end{equation*}
Similarly for the right-hand side, we have:
\begin{equation*}
\begin{split}
\prod_{i=1}^{f} {\left(1-\textbf{s}^{d_i}t^{n_i}\right)} &=
\prod_{\text{sgn}(n_i) = -1} {-\textbf{s}^{d_i}t^{n_i}}
 \ \prod_{i=1}^{f}{\left(1-\textbf{s}^{\text{sgn}(n_i)d_i}t^{|n_i|}\right)}.
\end{split}
\end{equation*}
Next, set \[\eta = \gamma \prod_{\text{sgn}(m_i) = -1}
{-\textbf{s}^{-c_i}t^{-m_i}}\prod_{\text{sgn}(n_i) = -1}
{-\textbf{s}^{d_i}t^{n_i}}\] so that our original equation may be
written as \[\prod_{i=1}^{e}
{\left(1-\textbf{s}^{\text{sgn}(m_i)c_i}t^{|m_i|}\right)} = \eta
\prod_{i=1}^{f} {\left(1-\textbf{s}^{\text{sgn}(n_i)
d_i}t^{|n_i|}\right)}.\]  Comparing the lowest degree term (with
respect to $t$) on both sides, it follows that $\eta = 1$. It is
enough, therefore, to prove the claim in the case when
\begin{equation}\label{polyequalityreduction}
\prod_{i=1}^{e} {\left(1-\textbf{s}^{c_i}t^{m_i}\right)} =
\prod_{i=1}^{f}{\left(1-\textbf{s}^{d_i}t^{n_i}\right)}
\end{equation}
and the $m_i,n_i$ are positive.  Without loss of generality,
suppose the lowest degree nonconstant term on both sides of
(\ref{polyequalityreduction}) is $t^{m_1}$ with coefficient
$-\textbf{s}^{c_1} - \cdots - \textbf{s}^{c_u}$ on the left and
$-\textbf{s}^{d_1} - \cdots - \textbf{s}^{d_v}$ on the right.
Here, $u$ (resp. $v$) corresponds to the number of $m_i$ (resp. $n_i$) with
$m_i = m_1$ (resp. $n_i = m_1$).

Since the set of distinct monomials $\{\textbf{s}^{c} : c \in C\}$
is a $\mathbb C$-vector space basis for the ring $\mathbb C[C]$, equality of the
$t^{m_1}$ coefficients above implies that $u = v$ and that up to
permutation, $\textbf{s}^{c_j} = \textbf{s}^{d_j}$ for $j =
1,\ldots,u$ (here is where we use that the characteristic of $\mathbb C$ is zero).  Lemma \ref{zerodivlemma} and induction complete the proof.
\end{proof}

\begin{lemma}\label{intmatrixlem}
Let $P = (p_{ij})$ be a $d$-by-$n$ integer matrix such that every
row has at least one nonzero integer.  Then, there exists
$\textbf{v} \in \mathbb Z^{n}$ such that the vector $P\textbf{v}$
does not contain a zero entry.
\end{lemma}

\begin{proof}
Let $P$ be a $d$-by-$n$ integer matrix as in the hypothesis of the
lemma, and for $h \in \mathbb Z$, let $\textbf{v}_h =
(1,h,h^2,\ldots,h^{n-1})^T$.  Assume, by way of contradiction,
that $P\textbf{v}$ contains a zero entry for all $\textbf{v} \in
\mathbb Z^{n}$.  Then, in particular, this is true for all
$\textbf{v}_h$ as above.  By the (infinite) pigeon-hole principle,
there exists an infinite set of $h \in \mathbb Z$ such that
(without loss of generality) the first entry of $P\textbf{v}_h$ is
zero.  But then, \[ f(h) := \sum\limits_{i = 1}^n {p_{1i} h^{i -
1} = 0}\] for infinitely many values of $h$.  It follows,
therefore, that $f(h)$ is the zero polynomial, contradicting our
hypothesis and completing the proof.
\end{proof}

Lemma \ref{intmatrixlem} will be useful in verifying the following
fact.

\begin{lemma}\label{homomorphismlemma}
Let $A$ be a finitely generated abelian group and $a_1,\ldots,a_d$
elements in $A$ of infinite order.  Then, there exists a
homomorphism $\phi:A \to \mathbb Z$ such that $\phi(a_i) \neq 0$
for all $i$.
\end{lemma}

\begin{proof}
Write $A = B \oplus C$, in which $C$ is a finite group and $B$ is
free of rank $n$.  If $n = 0$, then there are no elements of
infinite order; therefore, we may assume that the rank of $B$ is
positive.  Since $a_1,\ldots,a_d$ have infinite order, their
images in the natural projection $\pi: A \to B$ are nonzero.  It
follows that we may assume that $A$ is free and $a_i$ are nonzero
elements of $A$.

Let $e_1,\ldots,e_n$ be a basis for $A$, and write \[a_t =
p_{t1}e_1 + \cdots + p_{tn}e_n\] for (unique) integers $p_{ij} \in
\mathbb Z$. To determine a homomorphism $\phi: A \to \mathbb Z$ as
in the lemma, we must find integers $\phi(e_1),\ldots,\phi(e_n)$
such that
\begin{equation}
\begin{split}
 0 & \neq p_{11}\phi(e_1) + \cdots + p_{1n}\phi(e_n) \\
 & \ \cdots \cdots \cdots \cdots \cdots \cdots \cdots \cdots \\
 0 & \neq  p_{d1}\phi(e_1) + \cdots + p_{dn}\phi(e_n).
\end{split}
\end{equation}
This, of course, is precisely the consequence of Lemma
\ref{intmatrixlem} applied to the matrix $P = (p_{ij})$, finishing
the proof.
\end{proof}

Recall that a \textit{trivial unit} in the group ring $\mathbb C[A]$ is an
element of the form $\alpha \textbf{s}^a$ in which $0 \neq \alpha \in
\mathbb C$ and $a \in A$. The main content of Theorem
\ref{mainfactortheorem} is contained in the following result.  The
technique of embedding $\mathbb C[A]$ into a Laurent polynomial ring is
also used by Fried in \cite{fried}.

\begin{lemma}\label{factorlemmachar0}
Let $A$ be an abelian group. Two factorizations in $\mathbb C[A]$, \[
\prod_{i=1}^{e} \left(1-\textbf{s}^{g_i} \right) = \eta
\prod_{i=1}^{f} \left(1-\textbf{s}^{h_i}\right),\] in which $\eta$
is a trivial unit and $g_i,h_i \in A$ all have infinite order are
equal if and only if $e = f$ and there is some nonnegative integer
$p$ such that, up to permutation,
\begin{enumerate}
    \item $g_ i = h_i$ for $i = 1,\ldots,p$
    \item $g_i = -h_i$ for $i = p+1,\ldots,e$
    \item $\eta = (-1)^{e-p} \textbf{s}^{g_{p+1} + \cdots + g_e}$.
\end{enumerate}
\end{lemma}

\begin{proof}
The if-direction of the claim is a straightforward calculation.
Therefore, suppose that one has two factorizations as in the
lemma.  It is clear we may assume that $A$ is finitely generated.
By Lemma \ref{homomorphismlemma}, there exists a homomorphism
$\phi: A \to \mathbb Z$ such that $\phi(g_i),\phi(h_i) \neq 0$ for
all $i$.  The ring $\mathbb C[A]$ may be embedded into the Laurent ring,
$R = \mathbb C[A][t,t^{-1}],$ by way of
\[\psi \left(\sum_{i=1}^{m} { \alpha_i \textbf{s}^{a_i} } \right)
= \sum_{i=1}^{m} { \alpha_i \textbf{s}^{a_i}t^{\phi(a_i)} }.\]
Write $\eta = \alpha \textbf{s}^b$.  Then, applying this
homomorphism to the original factorization, we have
\[\prod_{i=1}^{e} \left(1-\textbf{s}^{g_i}t^{\phi(g_i)} \right) = \alpha
\textbf{s}^bt^{\phi(b)}\prod_{i=1}^{f}
\left(1-\textbf{s}^{h_i}t^{\phi(h_i)}\right).\]  Lemma
\ref{polyequality1var} now applies to give us that $e = f$ and
there is an integer $p$ such that up to permutation,
\begin{enumerate}
    \item $g_ i = h_i$ for $i = 1,\ldots,p$
    \item $g_i = -h_i$ for $i = p+1,\ldots,e$.
\end{enumerate}
We are therefore left with verifying statement (3) of the lemma.
Using Lemma \ref{zerodivlemma}, we may cancel equal terms in our
original factorization, leaving us with the following equation:
\begin{equation*}
\begin{split}
\prod_{i = p+1}^{e} (1-\textbf{s}^{g_i}) & =
\eta \prod_{i = p+1}^{e} (1-\textbf{s}^{-g_i}) \\
& = \eta(-1)^{e-p}\prod_{i = p+1}^{e} {\textbf{s}^{-g_i}}\prod_{i = p+1}^{e}{(1-\textbf{s}^{g_i})}. \\
\end{split}
\end{equation*}
Finally, one more application of Lemma \ref{zerodivlemma} gives us
that $\eta = (-1)^{e-p}\textbf{s}^{g_{p+1} +\cdots +g_e}$ as
desired. This finishes the proof.
\end{proof}

We may now prove Theorem \ref{mainfactortheorem}.

\begin{proof}[Proof of Theorem \ref{mainfactortheorem}]
Let \[ \textbf{s}^a \prod\limits_{i = 1}^e {\left(
{\textbf{s}^{u_i } - \textbf{s}^{v_i } } \right) = \alpha
\textbf{s}^b } \prod\limits_{i = 1}^f {\left( {\textbf{s}^{x_i } -
\textbf{s}^{y_i } } \right)}\] be two factorizations in the ring
$\mathbb C[Q]$. View this expression in $\mathbb C[A]$ and factor each element of
the form $\left( {\textbf{s}^{u } - \textbf{s}^{v } } \right)$ as
$\textbf{s}^{u }\left( {1 - \textbf{s}^{v-u } } \right)$.  By
assumption, each such $v-u$ has infinite order. Now, apply Lemma
\ref{factorlemmachar0}, giving us that $\alpha = \pm 1$, $e = f$,
and that after a permutation, for each $i$ either $\textbf{s}^{v_i
- u_i } = \textbf{s}^{y_i - x_i }$ or $\textbf{s}^{v_i - u_i } =
\textbf{s}^{ x_i - y_i}$. It easily follows from this that for
each $i$, there are elements $c_i,d_i \in Q$ such that
$\textbf{s}^{c_i}(\textbf{s}^{u_i}-\textbf{s}^{v_i}) = \pm
\textbf{s}^{d_i}(\textbf{s}^{x_i}-\textbf{s}^{y_i})$. This
completes the proof of the theorem.
\end{proof}

\section{Cyclic Resultants and Rational Functions}

We begin with some preliminaries concerning cyclic resultants. Let
$f(x) = a_0x^d + a_{1}x^{d-1} + \cdots + a_d$ be a degree $d$
polynomial over $\mathbb C$, and let the companion matrix for $f$ be given
by: \[A = \left[ {\begin{array}{*{20}c}
   0 & 0 &  \cdots  & 0 & { - a_d /a_0 }  \\
   1 & 0 &  \cdots  & 0 & { - a_{d - 1} /a_0 }  \\
   0 & 1 &  \cdots  & 0 & { - a_{d - 2} /a_0 }  \\
   0 &  \vdots  &  \ddots  &  \vdots  &  \vdots   \\
   0 & 0 &  \cdots  & 1 & { - a_1 /a_0 }  \\
\end{array} } \right].\]
Also, let $I$ denote the $d$-by-$d$ identity matrix.  Then, we may
write \cite[p. 77]{Cox}
\begin{equation}\label{rmdef}
r_m = a_0^m \text{det}\left( {A^m - I } \right).
\end{equation}
This equation can also be expressed as,
\begin{equation}\label{rmdef2}
r_m = a_0^m \prod\limits_{i = 1}^d {\left( {\alpha _i ^m -1}
\right)},
\end{equation}
in which $\alpha_1,\ldots,\alpha_d$ are the roots of $f(x)$.

Let $e_i(y_1,\ldots,y_d)$ be the $i$-th elementary symmetric
function in the variables $y_1,\ldots,y_d$ (we set $e_0 = 1$).
Then, we know that $a_i =
(-1)^{i}a_0e_i(\alpha_1,\ldots,\alpha_d)$ and that
\begin{equation}
r_m = a_0^m \sum_{i = 0}^d {(-1)^{i} e_{d - i} \left( {\alpha _1
^m , \ldots ,\alpha _d ^m } \right) }.
\end{equation}

We first record an auxiliary result.

\begin{lemma}\label{1elemma}
Let $F_k(z) = {\prod\limits_{1 \leq i_1 <  \cdots  < i_k \leq d}
{\left( {1 - a_0\alpha _{i_1 }  \cdots \alpha _{i_k }z } \right)}
}$ with $F_0(z) = 1- a_0z$.  Then,
\[\sum_{m = 1}^\infty {a_0^m e_k \left( {\alpha _1 ^m ,
\ldots ,\alpha _d ^m } \right)z^m } = -z \cdot \frac{{F_k'}}
{F_k},\] in which $F_k'$ denotes $\frac{{dF_k}} {{dz}}$.
\end{lemma}

\begin{proof}
For $k = 0$, the equation is easily verified.  When $k >0$, the
calculation is still fairly straightforward:

\begin{equation*}
\begin{split}
\sum_{m = 1}^\infty  {a_0^me_k \left( {\alpha_1^m , \ldots
,\alpha_d^m } \right)z^m } & = \sum_{m = 1}^\infty { \
\sum\limits_{i_1  <  \cdots  < i_k}
{a_0^m\alpha_{i_1 }^m  \cdots \alpha_{i_k }^m } \cdot z^m } \\
& = \sum\limits_{i_1  <  \cdots  < i_k} { \ \sum_{m = 1}^\infty
{a_0^m \alpha _{i_1 } ^m  \cdots \alpha_{i_k }^m \cdot z^m } } \\
& =  \sum\limits_{i_1  < \cdots < i_k} {\frac{a_0\alpha_{i_1 }
\cdots \alpha_{i_k }z} {{1 - a_0\alpha_{i_1 } \cdots \alpha_{i_k
}z }}} \\
& = \frac{{-z \cdot \frac{d} {{dz}}\left[ {\prod\limits_{i_1 <
\cdots  < i_k} {\left( {1 - a_0\alpha_{i_1 } \cdots \alpha_{i_k }z
} \right)} } \right]}} {{\prod\limits_{i_1 < \cdots  < i_k}
{\left( {1 - a_0\alpha_{i_1 } \cdots \alpha_{i_k }z } \right)} }} \\
& = -z \cdot \frac{{F_k'}} {F_k}.
\end{split}
\end{equation*}
\end{proof}

We are now ready to state and prove the rationality result mentioned in Section \ref{binfact}.

\begin{lemma}\label{rationalthm}
$R_f(z) = \sum\nolimits_{m = 1}^\infty  {r_m z^m }$ is a rational
function in $z$.
\end{lemma}

\begin{proof} We simply compute that
\begin{equation*}
\begin{split}
\sum_{m = 1}^\infty  {r_m z^m }  & = \sum_{m = 1}^\infty {\sum_{i
= 0}^d {(-1)^{i}a_0^me_{d - i} \left( {\alpha_1^m , \ldots
,\alpha_d^m } \right) }  \cdot z^m }
\\ & = \sum_{i = 0}^d {(-1)^{i} \sum_{m =
1}^\infty  {a_0^m e_{d - i} \left( {\alpha_1^m , \ldots
,\alpha_d^m } \right)}  \cdot z^m }  \\ & =  -z \cdot \sum_{i =
0}^d {(-1)^{i} \cdot \frac{{F_{d - i} '}} {{F_{d - i} }}}.
\end{split}
\end{equation*}
\end{proof}

Manipulating the expression for $R_f(z)$ occurring in Lemma
\ref{rationalthm}, we also have the following fact.

\begin{corollary}\label{gdcor}
If $d$ is even, let $G_d = \frac{{F_d F_{d - 2}  \cdots F_0 }}
{{F_{d - 1} F_{d - 3}  \cdots F_1 }}$ and if $d$ is odd, let $G_d
= \frac{{F_d F_{d - 2}  \cdots F_1 }} {{F_{d - 1} F_{d - 3} \cdots
F_0 }}$.  Then, \[\sum_{m = 1}^\infty  {r_m z^m }= -z
\frac{G_d'}{G_d}.\]
\end{corollary}

In particular, it follows that
\begin{equation}\label{expversion}
\exp{\left(- \sum_{m = 1}^\infty {r_m\frac{z^m}{m} } \right)}=
G_d.
\end{equation}

\begin{example}
Let $f(x) = x^2-5x+6 = (x-2)(x-3)$.  Then, $r_m = (2^m-1)(3^m-1)$
and $F_0(z) = 1-z$, $F_1(z) = (1-2z)(1-3z)$, $F_2(z) = 1-6z$.
Thus, \[R_f(z) = -z\left(\frac{F_2'}{F_2}- \frac{F_1'}{F_1} +
\frac{F_0'}{F_0}\right) = \frac{6z}{1-6z} -\frac{2z}{1-2z}-
\frac{3z}{1-3z} + \frac{z}{1-z}\] and \[\exp{\left(- \sum_{m =
1}^\infty {r_m\frac{z^m}{m} } \right)} =
\frac{(1-6z)(1-z)}{(1-2z)(1-3z)}.\qed\]
\end{example}

Following \cite{fried}, we discuss how to deal with absolute
values in the real case.  Let $f \in \mathbb R[x]$ have
degree $d$ such that the $r_m$ as defined above are all nonzero.
We examine the sign of $r_m$ using equation (\ref{rmdef2}). First
notice that a complex conjugate pair of roots of $f$ does not
affect the sign of $r_m$. A real root $\alpha$ of $f$ contributes
a sign factor of $+1$ if $\alpha > 1$, $-1$ if $-1 < \alpha < 1$,
and $(-1)^m$ if $\alpha < -1$.  Let $E$ be the number of zeroes of
$f$ in $(-1,1)$ and let $D$ be the number of zeroes in
$(-\infty,-1)$.  Also, set $\epsilon = (-1)^E$ and $\delta =
(-1)^D$.
 Then, it follows that \begin{equation}\label{rmabsvalues}
\frac{r_m}{|r_m|} = \epsilon \cdot \delta^{m}. \end{equation}  In particular,
\begin{equation}\label{rmrealcase}
|r_m| = \epsilon(\delta a_0)^m \prod\limits_{i = 1}^d {\left(
{\alpha _i ^m -1} \right)}.
\end{equation}
In other words, the sequence of $|r_m|$ is obtained by multiplying
each cyclic resultant of the polynomial $\tilde{f} := \delta f =
\delta a_0x^d + \delta a_1 x^{d-1}+\cdots+\delta a_d$ by
$\epsilon$. Denoting by $\widetilde{G}_d$ the rational function
determined by $\tilde{f}$ as in (\ref{gdcor}), it follows that
\begin{equation}\label{realexpformula}
\exp{\left(- \sum_{m = 1}^\infty {|r_m|\frac{z^m}{m} } \right)} =
\left(\widetilde{G}_d\right)^{\epsilon}.
\end{equation}

\section{Proofs of the Main Theorems}

Let $G$ be the multiplicative group generated by the roots
$\alpha_1,\ldots,\alpha_d$ of a polynomial $f$ for which $f(0) \neq 0$.  
We deal with the case when zero is a root of $f$ later.
Because of the multiplicative structure of $G$, 
we represent vector space basis elements of the group ring $\mathbb C[G]$ 
as $[\alpha]$, $\alpha
\in G$; multiplication is given by $[\alpha] \cdot [\beta] = [\alpha \beta]$.  
The divisor (in $\mathbb C[G]$) of the rational function $G_d$
defined by Corollary \ref{gdcor} is
\begin{equation}
( - 1)^{d+1} \left( {\sum\limits_{k\;\text{odd}}
{\;\sum\limits_{i_1 < \cdots  < i_k } {\left[ {\left( {a_0 \alpha
_{i_1 }  \cdots \alpha _{i_k } } \right)^{ - 1} } \right]} }  -
\sum\limits_{k\;\text{even}} {\;\sum\limits_{i_1  <  \cdots  < i_k
} {\left[ {\left( {a_0 \alpha _{i_1 }  \cdots \alpha _{i_k } }
\right)^{ - 1} } \right]} } } \right)
\end{equation}
\[= \left[ {a_0 ^{ - 1} } \right]\prod\limits_{i = 1}^d {\left(
{\left[ {\alpha _i ^{ - 1} } \right] - \left[ 1 \right]} \right)}.
\] Let us remark that for ease of presentation above, when $k = 0$, we
have assigned \[ {\sum\limits_{i_1  <  \cdots  < i_k } {\left[
{\left( {a_0 \alpha _{i_1 }  \cdots \alpha _{i_k } } \right)^{ -
1} } \right]} } = [a_0^{-1}],\] which corresponds to the factor of
$F_0(z) = 1-a_0z$ in $G_d$.  

Now, suppose that $f = x^{l}h$ in which $h(0) \neq 0$ and 
$h$ has degree $d$.  Then, from (\ref{rmdef2}), the
cyclic resultants of $f$ are given by $(-1)^l r_m(h)$.  
Examining equation (\ref{gdcor}) following 
Corollary \ref{gdcor}, it follows that the divisor of 
$G_d$ for $f$ is given by the divisor of the 
rational function
\[ \exp{\left(- \sum_{m = 1}^\infty {r_m(f)\frac{z^m}{m} } \right)} =
\left[ \exp{\left(- \sum_{m = 1}^\infty {r_m(h)\frac{z^m}{m} } \right) }  \right]^{(-1)^l}.\]
Let $\alpha_1,\ldots,\alpha_d$ be the roots of $h$.
By the discussion above, it therefore follows that the divisor of $G_d$ for $f$ is
\[ (-1)^l \left[ {a_0 ^{ - 1} } \right]\prod\limits_{i = 1}^d {\left(
{\left[ {\alpha _i ^{ - 1} } \right] - \left[ 1 \right]} \right)}.\]

With this computation in hand, we now
prove our main theorems.

\begin{proof}[Proof of Theorem \ref{characterizationthmchar0}]
Let $f$ and $g$ be polynomials as in the hypothesis, and suppose that 
the multiplicity of $0$ as a root of $f$ (resp. $g$) is $l_1$
(resp. $l_2$). Then, $f(x) = x^{l_1}(a_0x^{d_1}+\cdots+a_{d_1})$ and
$g(x) = x^{l_2}(b_0x^{d_2}+\cdots+b_{d_2})$ in which $a_0$ and
$b_0$ are not $0$. Let $\alpha_1,\ldots,\alpha_{d_1}$ and
$\beta_1,\ldots,\beta_{d_2}$ be the nonzero roots of $f$ and $g$,
respectively, and let $G$ be the multiplicative group generated by
these elements.  Since $f$ and $g$ both generate the same
sequence of cyclic resultants, it follows that the divisor (in the
group ring $\mathbb C[G]$) of their corresponding rational functions (see
(\ref{expversion})) are equal.  By above, such divisors factor,
giving us that \[(-1)^{d_1+l_1}[a_0^{-1}]\prod_{i=1}^{d_1}
\left([1]-[\alpha_i^{-1}] \right) =
(-1)^{d_2+l_2}[b_0^{-1}]\prod_{i=1}^{d_2} \left([1]-[\beta_i^{-1}]
\right).\]  Since we have assumed that $f$ and $g$ generate a set
of nonzero cyclic resultants, neither of them can have a root of
unity as a zero.  Therefore, Lemma \ref{factorlemmachar0} applies
to give us that $d := d_1 = d_2$ and that up to a permutation,
there is a nonnegative integer $p$ such that
\begin{enumerate}
    \item $\alpha_ i = \beta_i$ for $i = 1,\ldots,p$
    \item $\alpha_i = \beta_i^{-1}$ for $i = p+1,\ldots,d$
    \item $(-1)^{d-p} = (-1)^{l_2-l_1}$, $a_0b_0^{-1} = \beta_{p+1} \cdots
\beta_d$.
\end{enumerate}
Set $u(x) = (x-\beta_{p+1})\cdots(x-\beta_{d})$ which 
has $\deg(u) \equiv l_2-l_1 \ (\text{mod} \ 2)$, 
and let $v(x) = b_0(x-\beta_{1})\cdots(x-\beta_{p})$ (note
that if $p = 0$, then $v(x) = b_0$) so that $g(x) =
x^{l_2}v(x)u(x)$. Now,
\[u(x^{-1})x^{\text{deg}(u)} = (-1)^{d-p}\beta_{p+1} \cdots
\beta_d(x-\beta_{p+1}^{-1})\cdots(x-\beta_{d}^{-1}),\] and thus
\begin{equation*}
\begin{split}
f(x) & =
x^{l_1}a_0b_0^{-1}v(x)(x-\beta_{p+1}^{-1})\cdots(x-\beta_{d}^{-1})
\\
& = (-1)^{l_2-l_1}x^{l_1}v(x)u(x^{-1})x^{\text{deg}(u)}.\\
\end{split}
\end{equation*}
Finally, the converse is
straightforward from (\ref{rmdef2}), completing the proof
of the theorem.
\end{proof}

The proof of Theorem \ref{charthmrealcase} is similar, employing
equation (\ref{realexpformula}) in place of (\ref{expversion}).

\begin{proof}[Proof of Theorem \ref{charthmrealcase}]
Since multiplication of a real polynomial by a power of $x$ does
not change the absolute value of a cyclic resultant, we may assume
$f,g \in \mathbb R[x]$ have nonzero roots. The result now follows
from (\ref{realexpformula}) and the argument used to prove the
if-direction of Theorem \ref{characterizationthmchar0}.
\end{proof}

\section{Reconstructing dynamical systems from their zeta functions}\label{reconstruct}

In this section, we describe how to explicitly reconstruct a polynomial 
from its cyclic resultants.  For an ergodic toral endomorphism as 
in the introduction, sequences $|r_m|$ correspond to cardinalities 
of sets of periodic points.  In particular, the \textit{zeta function}, 
\[Z(T,z) = \exp{\left(- \sum_{m = 1}^\infty { |\text{Per}_m(T)|\frac{z^m}{m} } \right)} ,\] 
of the dynamical system 
in question is simply another way of writing equation (\ref{realexpformula}).

In many of the applications \cite{duistermaat, Zworski, kedlaya, 
knots}, the defining polynomial is reciprocal, and the techniques 
discussed here restrict easily to this special case.  Furthermore, 
since reciprocal polynomials are uniquely determined without 
any genericity assumptions (see Corollary \ref{friedtheorem} 
and Corollary \ref{recipcorollary}), the computational organization is simpler.

Let $f(x) = a_0x^d + a_{1}x^{d-1} + \cdots + a_d$ be a degree $d$
polynomial with indeterminate coefficients $a_i$.  We distinguish 
between two cases.  In the first situation, the variable $a_0$ is 
replaced by $1$ so that $f$ is monic; while in the second, 
we set $a_i = a_{d-i}$ for $i = 1,\ldots,d$ so that $f$ is reciprocal.  

Although the results mentioned in this paper only imply that the 
full sequence of cyclic resultants determine $f$ when it is 
(generic) monic or reciprocal, a finite number of resultants 
is sufficient.  Specifically, as detailed in forthcoming work 
\cite{hillarlevine}, it is shown that $2^{d+1}$ resultants are 
enough.  Empirical evidence suggests that this is far from 
tight, and a conjecture of Sturmfels and Zworski asserts the following.

\begin{conjecture}\label{sturmzworskiconj}
A generic monic polynomial $f(x) \in \mathbb C[x]$ of degree 
$d$ is determined by its 
first $d+1$ cyclic resultants.  Moreover, if $f$ is (non-monic) reciprocal 
of even degree $d$, then the number of resultants needed for inversion is given by $d/2 + 2$.
\end{conjecture}

A straightforward algorithm for inverting $N$ cyclic resultants 
is as follows.  Its correctness when $N = 2^{d+1}$ follows 
from \cite{Cox} and the results of \cite{hillarlevine}.

\begin{algorithm}\label{regalg}
(Specific reconstruction of a polynomial from its cyclic resultants) \\
Input:  Positive integer $d$ and a sequence of $r_1,\ldots,r_N \in \mathbb C$.  \\
Output:  The coefficients $a_i$ ($i=0,\ldots,d$) corresponding to $f$. 
\begin{enumerate}
\item  Compute a lexicographic Gr\"obner basis $\mathcal{G}$ 
for the ideal \[I = \langle r_1 - \text{Res}(f,x-1), \ldots, r_N - \text{Res}(f,x^N-1) \rangle.\]
\item Solve the resulting triangular system of equations for $a_i$ using back substitution.
\end{enumerate}
\qed
\end{algorithm}

If the data are given in terms of cyclic resultant absolute values 
(for the real case), then more care must be taken in implementing 
Algorithm \ref{regalg}.   Examining expression (\ref{rmabsvalues}), 
there are $2$ possible sequences of viable $r_m$ that come 
from a given sequence of (generically generated) cyclic 
resultant absolute values $|r_m|$; they are $\{|r_m|\}$ and 
$\{-|r_m|\}$.   By the uniqueness in Corollaries \ref{genericcormonic} 
and \ref{genericcorreal}, however, only one of these sequences 
can come from a monic polynomial.  Therefore, the corresponding 
modification is to run Algorithm \ref{regalg} on both these inputs.  
For one of these sequences, it will generate the Gr\"obner 
basis $\langle 1 \rangle$; while for the other, it will output the desired reconstruction.

Finding ``universal'' equations expressing the coefficients $a_i$ 
in terms of the resultants $r_i$ is also possible using a similar strategy.

\begin{algorithm}\label{universalalg}
(Formal reconstruction of a polynomial from its cyclic resultants) \\
Input:  Positive integers $d$ and $N$.  \\
Output:  Equations expressing $a_i$ ($i=0,\ldots,d$) parameterized by $r_1,\ldots,r_N$. 
\begin{enumerate}
\item Let $R = \mathbb Q[a_0,\ldots,a_d,r_1,\ldots,r_N]$ and let 
$\prec$ be any elimination term order with $\{a_i\}\prec \{r_j\}$.
\item  Compute the reduced Gr\"obner basis $\mathcal{G}$ for 
the ideal \[I = \langle r_1 - \text{Res}(f,x-1), \ldots, r_N - \text{Res}(f,x^N-1) \rangle.\]
\item Output a triangular system of equations for $a_i$ in terms of the $r_i$.
\end{enumerate}
\qed
\end{algorithm}

A few remarks concerning Algorithm \ref{universalalg} are in order. 
 If the $a_i$ are indeterminates, a monic polynomial with 
 coefficients $a_i$ will be generic.  Therefore, the first 
 $N = 2^{d+1}$ cyclic resultants of $f$ will determine it as a 
 polynomial in $x$ over an algebraic closure of 
 $\mathbb Q(a_1,\ldots,a_d)$.  It then follows from 
 general theory (for instance, quantifier elimination for 
 ACF, algebraically closed fields) that each $a_i$ can 
 be expressed as a rational function in the $r_i$ ($i = 1,\ldots,N$).  
 The same result holds for reciprocal polynomials with 
 indeterminate coefficients.  It is an interesting and difficult 
 problem to determine these rational functions for a given 
 $d$.  As motivation for future work on this problem, we 
 use Algorithm \ref{universalalg} to find these expressions 
 explicitly for several small cases.

When $f = a_0 x + a_1$ is linear, we need only two nonzero 
cyclic resultants to recover the coefficients $a_0, a_1$.  
An inversion is given by the formulae:  \[  a_0 = \frac{r_2^2 - r_1}{2r_1}, 
\  a_1 = \frac{-r_1^2 - r_2}{2r_1}.\]  In the quadratic case, a monic 
$f =  x^2 + a_1 x + a_2$ is also determined by two nonzero resultants:  
\begin{equation*}
\begin{split}
a_1 = \frac{r_1^2-r_2}{2r_1}, \ a_2 = \frac{r_1^2-2r_1+r_2}{2r_1}.  \\
\end{split}
\end{equation*}
When $f =  x^3 + a_1 x^2 + a_2 x + a_3$ has degree three, four 
resultants suffice, and inversion is given by: 
\begin{equation*}
\begin{split}
a_1 = &\  \frac{-12r_2r_1^3-12r_1r_2^2+3r_2^3-
r_2r_1^4-8r_2r_1r_3+6r_1^2r_4}{24r_2r_1^2}, \\ a_2 = 
& \ \frac{-r_1^2-2r_1+r_2}{2r_1}, \\  a_3 = & \ \frac{-3r_2^3+
r_2r_1^4+8r_2r_1r_3-6r_1^2r_4}{24r_1^2r_2}. \\
\end{split}
\end{equation*}
Reconstruction for $d=4$ is also possible using five resultants, 
however, the expressions are too cumbersome to list here.  

As a final example, we describe the reconstruction of a degree 
$6$ monic, reciprocal polynomial $f = x^6 + a_1x^5 + a_2x^4 + 
 a_3 x^3 +  a_2 x^2 + a_1x +1$ from its first four cyclic resultants:

\[P = -540\,{r_1}^{2}r_2\,r_4-13824\,{r_1}^{3}r_2+{ r_1}^{6}r_2+
27\,{r_2}^{3}{r_1}^{2}+9\,{r_1}^{ 4}{r_2}^{2}+27\,{r_2}^{4}-
432\,{r_1}^{3}{r_2}^{2}- \]\[648\,r_1\,{r_2}^{3}-72\,{r_1}^{5}r_2-
448\,r_3 \,{r_1}^{3}r_2+192\,r_3\,r_1\,{r_2}^{2}+108\,
{r_1}^{4}r_4+1536\,{r_1}^{2}r_2\,{\it
r_3}+\]\[2592\,{r_1}^{3}r_4+1728\,{r_1}^{4}r_2+5184\,{r_1}^{2} {r_2}^{2},\]

\[Q = {r_1}^{2} \left( -16\,r_3\,r_2+9\,r_4\,r_1
 \right),\]

\[R =  -648\,r_1\,{r_2}^{3}+27\,{r_2 }^{3}{r_1}^{2}+27\,{r_2}^{4}-
576\,r_3\,r_1\,{r_2}^{2}+2592\,{r_1}^{3}r_4+{r_1}^{6}r_2-
72\,{r_1}^{5}r_2+\] \[9\,{r_1}^{4}{r_2}^{2}+1728\,{r_1}^ {4}r_2-
432\,{r_1}^{3}{r_2}^{2}+ 320\,r_3\,{r_1
}^{3}r_2-324\,{r_1}^{4}r_4-13824\,{r_1}^{3}{\it r_2 }+\] 
\[5184\,{r_1}^{2}{r_2}^{2}+1536\,{r_1}^{2}r_2\,r_3-108\,{r_1}^{2}r_2\,r_4,\]

\[a_1 = {\frac {1}{192}}\,P/Q, \ a_2 = \,{\frac {-4\,r_1+{r_1}^{2}+
r_2}{4r_1}}, \ a_3 =\frac{-1}{96}R/Q.\]

\section{Acknowledgement}

We would like to thank Bernd Sturmfels and Maciej Zworski 
for bringing this problem
to our attention and for useful discussions.  We also thank 
the anonymous referees for 
helpful comments that improved exposition.



\begin{thebibliography}{1}

\bibitem{Cox}
D. Cox, J. Little, D. O'Shea, \emph{Using Algebraic Geometry},
Springer, New York, 1998.

\bibitem{duistermaat}
J.J. Duistermaat and V. Guillemin, \emph{The spectrum of positive
elliptic operators and periodic bicharacteristics}, Inv. Math. 25
(1975) 39-79.

\bibitem{Ward}
G. Everest and T. Ward. Heights of Polynomials and Entropy in
Algebraic Dynamics. Springer-Verlag London Ltd., London, 1999.

\bibitem{fried}
D. Fried, \emph{Cyclic resultants of reciprocal polynomials}, in
Holomorphic Dynamics (Mexico 1986), Lecture Notes in Math. 1345,
Springer Verlag, 1988, 124-128.

\bibitem{guillemin}
V. Guillemin, \emph{Wave trace invariants}, Duke Math. J. 83
(1996), 287-352.

\bibitem{hillarlevine}
C. Hillar and L. Levine, \emph{Polynomial recurrences and cyclic resultants}, submitted.

\bibitem{Zworski}
A. Iantchenko, J. Sj\"{o}strand, and M. Zworski, \emph{Birkhoff
normal forms in semi-classical inverse problems}, Math. Res. Lett. 9 (2002), 337-362.

\bibitem{kedlaya}
K. Kedlaya, \emph{Quantum computation of zeta functions of curves}, preprint.

\bibitem{CCA}
E. Miller and B. Sturmfels, \emph{Combinatorial Commutative
Algebra}, Springer, 2004.

\bibitem{purbhoo}
K. Purbhoo, \textit{A nullstellensatz for amoebas}, \\
http://math.berkeley.edu/\textasciitilde
kpurbhoo/papers/amoebas.pdf.

\bibitem{knots}
W. H. Stevens, \emph{Recursion formulas for some abelian knot
invariants}, Journal of Knot Theory and Its Ramifications, Vol. 9,
No. 3 (2000) 413-422.

\end{thebibliography}
\end{document}